\documentclass[11pt, a4paper]{amsart}

\usepackage{amsfonts,amsmath,amssymb, amscd}

\newtheorem{theorem}{Theorem}[section]
\newtheorem{lemma}[theorem]{Lemma}

\newtheorem{prop}[theorem]{Proposition}
\newtheorem{corollary}[theorem]{Corollary}

\theoremstyle{definition}
\newtheorem{rem}[theorem]{Remark}

\newcommand\pf{\begin{proof}}
\newcommand\epf{\end{proof}}

\newcommand{\cmdblackltimes}{\mathop{\raisebox{0.2ex}{\makebox[0.92em][l]{${\scriptstyle\blacktriangleright\mathrel{\mkern-4mu}<}$}}}}
\newcommand{\cmddotltimes}{\mathop{\raisebox{0.12ex}{$\shortmid$}\raisebox{0.2ex}{\makebox[0.86em][r]{${\scriptstyle\gtrdot\mathrel{\mkern-4mu}<}$}}}}

\numberwithin{equation}{section}

\title[Affineness of some quotient dur sheaves]
{Affineness of some quotient dur sheaves of a~super affine group}

\author{Akira Masuoka}
\address{Akira Masuoka: 
Institute of Mathematics, 
University of Tsukuba, 
Ibaraki 305-8571, Japan}
\email{akira@math.tsukuba.ac.jp}

\dedicatory{Dedicated to Professor Mitsuhiro Takeuchi in honor of his distinguished career}

\begin{document}

\begin{abstract}
We prove that given a super affine closed subgroup $H$ of a super affine group $G$
over a field $k$ of charctersitic $\mathrm{ch}~k \ne 2$, the dur $k$-sheaf $G\tilde{\tilde{/}} H$ of right cosets
is affine if the affine $k$-group $\overline{H}$ assocoiated to $H$ is (a) reductive or (b) pro-finite.
Especially when $G$ is algebraic, the result in Case (a) gives rise to a positive answer to Brundan's
question which was recently discussed by Zubkov \cite{Z}.
\end{abstract}

\maketitle



\hspace{3cm}

\section{Introduction}\label{Introduction}

Throughout this paper we work over a fixed field $k$, whose characteristic $\mathrm{ch}~k$
is supposed to differ from 2 unless otherwise stated. 
Let $\mathbf{S}$ denote the category of vector spaces graded by the group 
$\mathbb{Z}_2 =\{0, 1\}$. 
It is a $k$-linear tensor category given the canonical symmetry
\begin{equation}\label{(1.1)}
V \otimes W \overset{\simeq}{\longrightarrow} W \otimes V,
\ v \otimes w \mapsto (-1)^{|v||w|}w \otimes v, 
\end{equation}
where $v(\in V), w(\in W)$ are supposed to be homogeneous elements of degree $|v|, |w|$.
Algebraic systems in the vector space category are generalized to those in $\mathbf{S}$.
The generalized latter are called with \emph{super} prefixed; 
for example, commutative (Hopf) algebras in $\mathbf{S}$ are called 
\emph{super commutative (Hopf) algebras}. 
A \emph{super $k$-functor} (resp., \emph{super $k$-group functor}) is a set-valued (resp., 
group-valued) functor on the category of super commutative algebras. 
Generalizing the subject worked out by Demazure and Gabriel \cite{DG}, Zubkov \cite{Z} 
recently defined the notion of (dur) $k$-sheaves in the super context, and successfully 
associated a (dur) $k$-sheaf with universal property to every super $k$-functor.

A super $k$-(group) functor is said to be \emph{affine} if it is representable. 
It is then represented by a super commutative (Hopf) algebra, say $A$, and the functor 
is denoted by $\mathrm{Sp}~A$ (denoted by $\mathrm{SSp}~A$ in \cite{Z}). 
An affine super $k$-group functor is called a \emph{super affine $k$-group}.
Let $G = \mathrm{Sp}~A$ be as such, where $A = A_0 \bigoplus A_1$ is a super commutative 
Hopf algebra.
It is said to be \emph{algebraic} (resp., \emph{finite}) if $A$ is finitely generated 
as an algebra (resp., finite-dimensional).
This $G$, restricted to the category of ordinary commutative algebras, gives rise to the
affine $k$-group $\overline{G} = \mathrm{Sp}~\overline{A}$ which is represented by the 
largest quotient ordinary Hopf algebra
\begin{equation}\label{(1.2)}
\overline{A} := A_0/A_1^{2}
\end{equation}
of $A$; see \cite[p.298]{M}. We call $\overline{G}$ the affine $k$-group \emph{associated to} $G$.

Our main results of this paper are the following two.

\begin{theorem}\label{1.1}
Let $G$ be a super affine $k$-group, and let $H$ is a super affine closed subgroup of $G$.
Then the dur $k$-sheaf $G\tilde{\tilde{/}} H$ associated to the super $k$-functor
$G/H$ of right cosets is affine, if the affine $k$-group $\overline{H}$ associated to $H$
is (a)~reductive or (b)~pro-finite.
\end{theorem}

\begin{corollary}\label{1.2}
Let $G, H$ be as above, and assume that $G$ is algebraic.
Then the super $k$-sheaf $G\tilde{/} H$ associated to $G/H$ is affine, and coincides with
$G\tilde{\tilde{/}} H$, if $\overline{H}$ is (a)~reductive or (b)~finite.
\end{corollary}

The last result in Case (a) answers in the positive Brundan's question discussed by Zubkov
\cite{Z}, who proved the same result as Corollary \ref{1.2} in the restricted situation
when (a) $\mathrm{ch}~k > 2$  and $\overline{H}$ is reductive or (b) $G$ is finite.

If $G$ is algebraic in the situation above, then $H$ and hence $\overline{H}$ are algebraic,
in which case Condition (b) in Theorem \ref{1.1}, which is restated so as by (b) in 
Theorem \ref{1.4} below, is equivalent to (b) in Corollary \ref{1.2}. 
Therefore, the corollary follows from Theorem \ref{1.1} and the following.

\begin{prop}\label{1.3}
Let $G, H$ be as in Theorem \ref{1.1}. Assume that $G$ is algebraic. If $G\tilde{\tilde{/}} H$
is affine, then $G\tilde{/} H$ is affine, and coincides with $G\tilde{\tilde{/}} H$.
\end{prop}

Our three results above remain true if $G/H$ is replaced by the super $k$-functor $H \backslash G$
of left cosets; 
to see this, apply each of the results to the opposite group $G^{\mathrm{op}}$.

Let $G = \mathrm{Sp}~A$ be a super affine $k$-group. Every super affine closed subgroup $H$
of $G$ arises uniquely from a quotient super Hopf algebra $A \to D$ of $A$ so that 
$H = \mathrm{Sp}~D$.
As a genralization of Takeuchi's Theorem \cite[Theorem 10]{T2}, it is proved by Zubkov 
\cite[Theorem 5.2]{Z} that $G\tilde{\tilde{/}} H$ (or equivalently, $G\tilde{\tilde{\backslash}} H$)
is affine if and only if $A$ is faithfully coflat, regarded as a right or equivalently, left
$D$-comodule along the quotient map $A \to D$. 
If this last condtion is satisfied, we say that
\begin{equation}\label{(1.3)}
A \to D\ \emph{is\ faithfully\ coflat};
\end{equation}
see Proposition \ref{2.2} below for some equivalent conditions. 
In virtue of the equivalence stated above, Theorem \ref{1.1} is translated into Hopf-algebra
language as follows.

\begin{theorem}\label{1.4}
A quotient $A \to D$ of a super commutative Hopf algebra $A$ is faithfully coflat if the commutative
Hopf algebra $\overline{D} = D_0/D_1^{2}$ is (a)~cosemisimple or (b)~a directed union of 
finite-dimensional Hopf subalgebras.
\end{theorem}

Generalizing \cite[III, Sect.3, 7.2]{DG}, the main theorem, Theorem 6.2, of \cite{Z} states that
if a super affine closed subgroup $H$ of a super affine $k$-group $G$ is normal, then 
$G\tilde{\tilde{/}} H$ is affine and is a super affine $k$-group. 
A Hopf-algebraic counterpart of this result as well as of some others from \cite{Z} (published
2009) had been proved by the author in the article \cite{M} published 2005.
In the article \cite{M} just cited, an important role was played by the Tensor Product Decomposition
Theorem, which will be reproduced in Section \ref{4} as Theorem \ref{4.1}, and whose proof will be
given there in a refined form because a part of the original proof was not quite well.
That theorem plays an important role in this paper as well, to prove Theorem \ref{1.4} above.
The proof of Theorem \ref{1.4} is given in the last Section \ref{5}, while Proposition \ref{1.3}
is proved in Section \ref{3}. 
The two results imply the remaining Theorem \ref{1.1} and Corollary \ref{1.2}, as was already
noted.

\section{Some basic results on super (co)algebras}\label{2}

As in \cite{M}, we will often write $\mathbb{Z}_2$ to denote the group (Hopf) algebra $k\mathbb{Z}_2$.

Given an ordinary (resp., super) algebra $R$, we let ${}_{R}\mathbf{M}$, $\mathbf{M}_R$ (resp.,
${}_{R}\mathbf{S}$, $\mathbf{S}_R$) denote the categories of left and respectively, right $R$-modules
(resp., those modules in $\mathbf{S}$). 
If $R$ is super, it is naturally regarded as a right $\mathbb{Z}_2$-module algebra, which constitutes 
the algebra $\widetilde{R} = \mathbb{Z}_2 \ltimes R$ of smash product.
Note that ${}_{R}\mathbf{S} = {}_{\widetilde{R}}\mathbf{M}$, $\mathbf{S}_R = \mathbf{M}_{\widetilde{R}}$.

\begin{prop}\label{2.1}
For a super algebra $R$, the following are equivalent:
\begin{itemize}

\item[(a)] 
$R$ is a right Noetherian ring;
 
\item[(b)] 
$R$ is a Noetherian object in $\mathbf{S}_R$, or in other words, super right ideals in $R$ satisfy
the descending chain condition;

\item[(c)] 
$\widetilde{R} = \mathbb{Z}_2 \ltimes R$ is a right Noetherian ring.
\end{itemize}
The parallel conditions with `right' replaced by `left' are equivalent to each other.
\end{prop}

\begin{proof}
(a) $\Rightarrow$ (b): Obvious.

(b) $\Rightarrow$ (c): Let $e_i~(i = 1, 2)$ denote the primitive idempotents in the
group (Hopf) algebra $k\mathbb{Z}_2$. 
Suppose that $e_1$ corresponds to the counit $\varepsilon$, so that 
$\varepsilon(e_i) = \delta_{i,1}.$
We see that $\widetilde{R} = \bigoplus_{i=1,2}e_{i} \otimes R$ in $\mathbf{S}_R$, 
$e_{1} \otimes R \simeq R$ in $\mathbf{S}_R$, and $e_{2} \otimes R$ is a degree shift of $R$.
Therefore, (b) implies that $e_{i} \otimes R~(i = 1, 2)$  are both Noetherian in $\mathbf{S}_R$,
which in turn implies (c).

(c) $\Rightarrow$ (a): This is seen if one introduces to $\widetilde{R}$, just as in 
\cite[Page 304, lines 17--18]{M}, an alternative, but isomorphic
structure in $\mathbf{M}_{\widetilde{R}}$.

We have worked above in the `right case'. In the `left' case, discuss in parallel, using a mirror.
\end{proof}

Lemma 5.1 of \cite{M} characterizes faithfully flat or projective modules over a super
algebra. 
We will dualize the result as far as will be needed. 
Given an ordinary (resp., super) coalgebra $C$, we let ${}^{C}\mathbf{M}$, $\mathbf{M}^{C}$ (resp.,
${}^{C}\mathbf{S}$, $\mathbf{S}^{C}$) denote the categories of left and respectively, right $C$-comodules
(resp., those comodules in $\mathbf{S}$). 
If $C$ is super, it constitutes the coalgebra $\widetilde{C} = \mathbb{Z}_2 \cmdblackltimes C$ of smash coproduct.
Note that ${}^{C}\mathbf{S} = {}^{\widetilde{C}}\mathbf{M}$, $\mathbf{S}^{C} = \mathbf{M}^{\widetilde{C}}$.

\begin{prop}\label{2.2}
For an object $V$ in $\mathbf{S}^{C}$, the following are equivalent:
\begin{itemize}

\item[(a)]
$V$ is (an) injective (cogenerator) in $\mathbf{M}^{C}$;

\item[(b)]
$V$ is (faithfully) coflat in $\mathbf{M}^{C}$;

\item[(c)]
The cotensor product functor $V \square_{C} : {}^{C}\mathbf{S} \to \mathbf{S}$ is (faithfully) exact.
\end{itemize}
A parallel result holds true for every object in ${}^{C}\mathbf{S}$.
\end{prop}

\begin{proof}
(a) $\Leftrightarrow$ (b): This is due to Takeuchi \cite[Proposition~A.2.1]{T1}.

(b) $\Leftrightarrow$ (c): Dualize the proof of \cite[Lemma~5.1(1)]{M}. 
In fact, (b) $\Rightarrow$ (c) is obvious. 
For the converse, let $\widetilde{V} = \mathbb{Z}_2 \otimes V$ denote the tensor 
product in $\mathbf{S}$;~this is an object in $\mathbf{S}^{C}$, given the right C-comodule
structure arising from that on $V$. 
If $W \in~^{C}\mathbf{S}$, we see that $\varepsilon \otimes \mathrm{id} \otimes \mathrm{id} :
\mathbb{Z}_2 \otimes V \otimes W \to V \otimes W$ induces a natural isomorphism 
$\widetilde{V} \square_{\widetilde{C}} W \overset{\simeq}{\longrightarrow} V \square_{C} W$.  
Therefore, (c) implies that $\widetilde{V}$ is (faithfully) coflat in $\mathbf{M}^{\widetilde{C}}$, 
which in turn implies (b), since $\widetilde{C}$ is faithfully coflat in $\mathbf{M}^{C}$.
\end{proof}

Next, we give a supplementary result to \cite{M}; see \cite{T3} for a more substantial, supplementary
result to \cite{M}, which concerns structure of super cocommutative Hopf algebras. 
Let $R$ be a super algebra.  Let $R^{\circ}$ denote the dual super coalgebra of $R$ as defined
in \cite[p.290, line~16]{M}; this consists of those elements in the dual vector space $R^*$ 
which annihilate some super ideal in $R$ of cofinite dimension.

\begin{lemma}\label{2.3}
This $R^{\circ}$ coincides with the dual coalgebra of $R$ as defined in \cite[p.109]{Sw} in
non-super context.

\end{lemma}

\begin{proof}
The latter is defined to be the pullback of $R^* \otimes R^*$ along the dual
$R^* \to (R \otimes R)^*$ of the product map, which is therefore $\mathbb{Z}_2$-graded.
Each homogeneous element in it annihilates some super ideal in $R$ of cofinite dimension, as
is seen from the proof of \cite[Proposition~6.0.3]{Sw}, 3) $\Rightarrow$ 4) $\Rightarrow$ 1).
This proves the lemma.
\end{proof}

\begin{rem}\label{2.4}
The three results above are generalized as follows, with $\mathbb{Z}_2$ replaced by more 
general Hopf algebras.  
Here the characteristic $\mathrm{ch}~k$ may be arbitrary.

(1) Let $J$ be a finite-dimensional semisimple commutative Hopf algebra, and let $R$ be a 
right, say, $J$-module algebra, which constitutes the algebra $J \ltimes R$ of smash product. 
Proposition \ref{2.1} holds for this $R$, with $\widetilde{R}, \mathbf{S}_R$ replaced by 
$J \ltimes R, \mathbf{M}_{J \ltimes R}$.

(2) Let $J$ be a Hopf algebra with bijective antipode, and let $C$ be a right, say, 
$J$-comodule coalgebra, which constitutes the coalgebra  $J \cmdblackltimes C$ of smash coproduct.  
Proposition \ref{2.2} holds in the genralized situation that  ${}^{C}\mathbf{S}, \mathbf{S}^C$ 
are replaced by the categories ${}^C(\mathbf{M}^{J}),~(\mathbf{M}^{J})^C (= \mathbf{M}^{J \cmdblackltimes C})$
of left and respectively, right $C$-comodules in the tensor category $\mathbf{M}^{J}$.

(3) Let $R$ be an algebra graded by a finite group $G$. 
Lemma \ref{2.3} is generalized so that if an element in $R^*$ annihilates some ideal in $R$ 
of cofinite dimension, it necessarily annihilates some $G$-graded ideal of cofinite dimension.
\end{rem}

\section{Proof of Proposition \ref{1.3}}\label{3}

A super commutative algebra $R$ is said to be \emph{Noetherian} if super ideals in $R$ satisfy 
the descending chain condition. 
By Proposition \ref{2.1}, this last condition is equivalent to that $R$ is left or equivalently, 
right Noetherian ring.

\begin{lemma}\label{3.1}
Let $R \to A$ be a map of super commutative algebras, with which $A$ is regarded as 
a super algebra over $R$.  
Assume that $A$ is finitely generated over $R$, and $R$ is Noetherian.  
Then, $A$ is Noetherian, and is finitely presented over $R$ in the super sense as defined in 
\cite[Page 721, line 4--6]{Z}.
\end{lemma}

\begin{proof}
Let $P = k[x_1,..., x_n]$ denote a polynomial algebra in finite indeterminates, and 
let $T = \wedge(V)$ denote the exteior algebra of a finite-dimensional vector space $V$. 
It suffices to prove that the super commutative algebra $R \otimes P \otimes T$ is 
Noetherian, since $A$ is a homomorphic image of such a super algebra over $R$.
By (the proof of) Hilbert's Basis Theorem, $R \otimes P$ is Noetherian. 
This implies each sub-quotient $R \otimes P \otimes \wedge^{i}(V)$ of $R \otimes P \otimes T$        
is Noetherian, which in turn implies the desired Noetherian property.
\end{proof}

Let $A$ be a super commutative Hopf algebra, whose coalgebra structure maps will be denoted by
\begin{equation}\label{(3.1)}
\Delta : A \to A \otimes A,\ \Delta (a) = \sum a_1 \otimes a_2;\ \varepsilon : A \to k. 
\end{equation}
A super left or right coideal subalgebra $B \subset A$ is said to be \emph{faithfully flat} 
if $A$ is faithfully flat as a left or equivalently, right $B$-module; see 
\cite[Corollary 5.5]{M} also for further equivalent conditions. 
Recall from \cite[Proposition 5.6]{M} the following.

\begin{prop}\label{3.2}
Fix a super commutative Hopf algebra $A$. 
Then the faithfully flat super left coideal subalgebras $B \subset A$ and the faithfully coflat 
quotient super Hopf algebras $A \to D$ (see \eqref{(1.3)}) are in one-to-one correspondence, by
\begin{equation}\label{(3.2)}
B \mapsto A/B^{+}A,\ D \mapsto A^{\mathrm{co}D}. 
\end{equation}
\end{prop}

Here, $B^{+} = B \cap \mathrm{Ker}~\varepsilon$. Given a quotient $A \to D,~a \mapsto \overline{a}$, 
$A^{\mathrm{co}D}$ is defined by
\begin{equation}\label{(3.3)}
A^{\mathrm{co}D} = \{ a \in A \mid \sum a_1 \otimes \overline{a}_2 = a \otimes \Bar{1} \},
\end{equation}
the right $D$-coinvariants in $A$.

\subsection*{Proof of Proposition 1.3.}
Let $A \to D$ be a faithfully coflat quotient of a super commutative Hopf algebra $A$. 
This gives rise to a closed embedding $H := \mathrm{Sp}~D \hookrightarrow G := \mathrm{Sp}~A$
of super affine $k$-groups with $G\tilde{\tilde{/}} H$ affine; see \cite[Theorem~5.2]{Z}.  
Set $B = A^{\mathrm{co}D}$. 
By Proposition \ref{3.2}, $B \subset A$ is faithfully flat and $A/B^{+}A = D$.  
It follows by \cite[Proposition 1.3]{M} that $N \mapsto N \otimes_{B} A$ gives a category equivalence
\begin{equation}\label{(3.4)}
\mathbf{S}_{B} \overset{\approx}{\longrightarrow} \mathbf{S}_{A}^{D},
\end{equation}
where $\mathbf{S}_{A}^{D}$ denotes the category of ($D, A$)-Hopf modules in $\mathbf{S}$.

Assume that $G$ is algebraic, or $A$ is finitely generated. 
Since $A$ is then Noetherian by Lemma \ref{3.1}, it is a Noetherian object in $\mathbf{S}_{A}^{D}$, 
whence the corresponding $B$ in $\mathbf{S}_{B}$ is Noetherian.  
This, combined with Lemma \ref{3.1} and \cite[Proposition 5.2]{Z}, implies the desired result.\ $\square$

\section{Tensor product decomposition theorem, revisited}\label{4}

Throughout this section we let $A = A_0 \oplus A_1$ denote a super commutative Hopf algebra. 
Recall from \eqref{(1.2)} the definition of the quotient ordinary Hopf algebra $\overline{A}$
of $A$.  
As in \cite{M}, we define
\begin{equation}\label{(4.1)}
W^A := A_{1}/A_{0}^{+}A_{1},
\end{equation} 
where $A_{0}^{+} = A_0 \cap \mathrm{Ker}~\varepsilon$; this is the odd part of the cotangent
space of the super affine $k$-group $\mathrm{Sp}~A$ at unity.  
The assignments $A \mapsto \overline{A},\ A \mapsto W^A$ are functorial.  
Recall that for any vector spavce $V$, the exterior algebra $\wedge(V)$ forms a super commutative
and cocommutative Hopf algebra in which every element in $V$ is an odd primitive.  
We reproduce from \cite{M} the following.

\begin{theorem}\label{4.1}
{\rm (\cite[Theorem 4.5, Remark 4.8]{M})}

(1) There is a counit-preserving isomorphism 
$A \overset{\simeq}{\longrightarrow} \overline{A} \otimes \wedge(W^A)$ of super 
left $\overline{A}$-comodule algebras.

(2) Given a map $f : A \to D$ of super commutative Hopf algebras, isomorphisms 
$A \overset{\simeq}{\longrightarrow} \overline{A} \otimes \wedge(W^A)$, 
$D \overset{\simeq}{\longrightarrow} \overline{D} \otimes \wedge(W^D)$ such as 
claimed above can be chosen so that they make the following diagram commute:
\begin{equation}\label{(4.2)}
\begin{CD}
A @>{\simeq}>> \overline{A} \otimes \wedge(W^A) \\
@V{f}VV @VV{\overline{f} \otimes \wedge(W^f)}V \\
D @>{\simeq}>> \overline{D} \otimes \wedge(W^D)
\end{CD}
\end{equation}
\end{theorem}

The original proof given in \cite{M} was divided into 0-affine case and general case.
As the proof in the latter case was not quite well, we will refine it below. 
To start we need to recall some argument in 0-affine case.

The dual super coalgebra $A^{\circ}$ of $A$ (see Lemma \ref{2.3}) is a super cocommutative
Hopf algebra. Let
\begin{equation}\label{(4.3)}
V_{A^{\circ}} := P(A^{\circ})_1
\end{equation}
denote the vector space of all odd primitives in $A^{\circ}$.

Assume that $A$ is 0-\emph{affine} \cite[Definition 4.2]{M} in the sense that
$A$ is finitely generated over the central subalgebra $A_0$. This implies 
$\mathrm{dim}~W^A < \infty$; see \cite[Proposition 4.4]{M}. 
In this case, $V_{A^{\circ}}$ and $W^A$ are naturally dual to each other. 
Moreover, the finite-dimensional super Hopf algebras $\wedge(V_{A^{\circ}})$ and 
$\wedge(W^A)$ are dual to each other; see \cite[Proposition 4.3(2) and Page~291, 
lines~1--3]{M}. 
Choose arbitrarily a totally ordered basis $X$ of $V_{A^{\circ}}$, and let $\iota_{X} 
: \wedge(V_{A^{\circ}}) \to A^{\circ}$ denote the unit-preserving super coalgebra map
defined by
\begin{equation}\label{(4.4)}
\iota_{X}(x_{\lambda} \wedge x_{\mu} \wedge ... \wedge x_{\nu}) = x_{\lambda}x_{\mu}~...~x_{\nu},
\end{equation}
where $x_{\lambda} < x_{\mu} <...< x_{\nu}$ in $X$. Let $\rho_{X} : A \to \wedge(W^A)$ 
denote the composite
\begin{equation}\label{(4.5)}
A \longrightarrow (A^{\circ})^* \overset{{\iota_{X}}^{*}}{\longrightarrow} \wedge(V_{A^{\circ}})^{*}
\overset{\simeq}{\longrightarrow} \wedge(W^A),
\end{equation}
in which the first arrow is the canonical map. 
It is proved in \cite[Proof in 0-affine case, pp. 300--301]{M} that
\begin{equation}\label{(4.6)}
\psi_X : A \to \overline{A} \otimes \wedge(W^A),\ \psi_X(a) = \sum \overline{a}_1 \otimes \rho_{X}(a_2)
\end{equation}
is necessarily such an isomorphism as claimed by Part 1 of the theorem above; 
as in \cite{M}, an isomorphism of this form will be said to be \emph{admissible}. 
We will identify $X$  with the dual basis of $W^A$ given the corresponding total order.

Given a map $f : A \to D$ to another 0-affine super commutative Hopf algebra $D$, 
choose totally ordered bases $X$ of $W^A$, and $Y$ of $W^D$ so that $W^f : W^A \to W^D$
restricts to a map $X \to Y \cup \{ 0 \}$ which strictly preserves the order on 
$X \setminus (W^f)^{-1}(0)$. 
(Note that the dual map $V_{f^{\circ}} : V_{D^{\circ}} \to V_{A^{\circ}}$ has the same property.) 
One then sees that the admissible isomorphisms $\psi_X,~\psi_Y$ make the diagram \eqref{(4.2)}
commute. 
It follows that if $f$ is an inclusion, in particular, then $\overline{f} : \overline{A} \to 
\overline{D},\ W^f : W^A \to W^D$ are both injections, through which we will regard so as
\begin{equation}\label{(4.7)}
\overline{A} \subset \overline{D},\ \ W^A \subset W^D,\ X \subset Y~(\mathrm{as~ordered~sets}).
\end{equation} 
In this case we write
\begin{equation}\label{(4.8)}
(A,~X) \subset (D,~Y).
\end{equation} 

\subsection*{Proof in general case.}
Let $A$ be in general.

(1) Let $\mathcal{F} = \mathcal{F}_A$ denote the set of those pairs $(B,~X)$ in which $B$ 
is a 0-affine super Hopf subalgebra of $A$, and $X$ is a totally ordered basis of $W^B$. 
This set is ordered with respect to the $\subset$ defined in \eqref{(4.8)}.  
Given a directed subset $\{ (B_{\alpha}, X_{\alpha}) \}_{\alpha}$ of $\mathcal{F}$, 
the directed union $B := \bigcup_{\alpha} B_{\alpha}$ is a super Hopf subalgebra of $A$ 
such that $\overline{B} = \bigcup_{\alpha} \overline{B}_{\alpha},\ W^B = \bigcup_{\alpha} W^{B_{\alpha}}$; 
see \eqref{(4.7)} and \cite[Proposition 4.3(3)]{M}.  
This implies that the inductive limit $\underrightarrow{\mathrm{lim}}_{\alpha}~\psi_{X_{\alpha}}$ gives 
an isomorphism on $B$ such as claimed by the theorem.

By Zorn's Lemma we have a directed subset $\mathcal{G} = \{ (B_{\alpha}, X_{\alpha}) \}_{\alpha}$ of 
$\mathcal{F}$ which is maximal with respect to inclusion. 
Set $B = \bigcup_{\alpha} B_{\alpha}$. 
The argument in the preceding paragraph shows that it suffices to prove $B = A$. 
Assume $B \subsetneq A$ on the contrary. 
As in the original proof we have a pair $(B',~X')$ in $\mathcal{F}$ and a 0-affine super Hopf subalgebra
$C \subset A$ such that $B' = B \cap C,\ C \nsubseteq B$. 
Extend the basis $X'$ of $W^{B'}$ to a basis $X' \sqcup Y$ (disjoin union) of $W^C$. 
For each $\alpha$ such that $(B_{\alpha},~X_{\alpha}) \supset (B',~X')$, consider the pair
\begin{equation}\label{(4.9)}
(B_{\alpha}C,~X_{\alpha} \sqcup Y).
\end{equation}
This $B_{\alpha}C$ is a 0-affine super Hopf subalgebra of $A$. 
Since $B_{\alpha}~\cap~C = B'$, it follows from the proof of \cite[Proposition 4.7]{M} that 
$W^{B_{\alpha}C} = W^{B_{\alpha}} +~W^C$, $W^{B_{\alpha}} \cap~W^C = W^{B'}$.  
Therefore, the $X_{\alpha} \sqcup Y$ in \eqref{(4.9)} is a basis of $W^{B_{\alpha}C}$.  
Choose arbitrarily a total order on $Y$, and extend the orders on $X_{\alpha}, Y$ to 
a total order on $X_{\alpha} \sqcup Y$, uniquely so that $X_{\alpha} < Y$ (that is, 
$x < y$ if $x \in X_{\alpha}, y \in Y$). 
Then the pair is in $\mathcal{F}$. 
The subset $\mathcal{G}$ joined with all pairs given in \eqref{(4.9)} forms a directed 
subset of $\mathcal{F}$ which properly includes, contradicting the maximality of $\mathcal{G}$. 
Therefore, we must have $B = A$, as desired.

(2) Let $f : A \to D$ be as in the theorem. Let $\mathcal{F}'_{A}$ denote the subset of $\mathcal{F}_{A}$ 
consisting of those  pairs $(B,~X)$ in which $X$ decomposes as $X = X' \sqcup X''$ so that 
$W^f : W^A \to W^D$ is injective on the subspace spanned by $X'$, and vanishes on $X''$. 
Let $\{ (B_{\alpha}, X_{\alpha}) \}_{\alpha}$ be a maximal directed subset of $\mathcal{F}'_{A}$, 
which exists by Zorn's Lemma. 
By modifying the proof of (1) above, it follows that $A = \bigcup_{\alpha} B_{\alpha}$. 
We see that $\mathcal{E} := \{ (f(B_{\alpha}),~f(X'_{\alpha})) \}_{\alpha}$ is a directed subset 
of $\mathcal{F}_{D}$, where we suppose $X_{\alpha} = X'_{\alpha} \sqcup X''_{\alpha}$ as above, 
and $f(X'_{\alpha})$ has the total order inherited from $X'_{\alpha}$. 
Moreover, the isomorphisms $\underrightarrow{\mathrm{lim}}_{\alpha}~\psi_{X_{\alpha}}$ on $A$, and
$\underrightarrow{\mathrm{lim}}_{\alpha}~\psi_{f(X'_{\alpha})}$ on $f(A) = \bigcup_{\alpha} f(B_{\alpha})$ 
make the diagram \eqref{(4.2)} with $D$ replaced by $f(A)$ commute. 
Again as above, we can choose a maximal directed subset $\{ (E_{\beta},~Y_{\beta}) \}_{\beta}$ of 
$\mathcal{F}_{D}$ which includes $\mathcal{E}$, and then necessarily have $\bigcup_{\beta} E_{\beta}= D$. 
The proof completes if we replace the last isomorphism on $f(A)$ with the isomorphism 
$\underrightarrow{\mathrm{lim}}_{\beta}~\psi_{Y_{\beta}}$ on $D$.\ $\square$

\section{Proof of Theorem \ref{1.4}}\label{5}

\begin{lemma}\label{5.1}
Let $A \to D$ be a quotient of a super commutative Hopf algebra $A$. 
Assume that $A = \bigcup_{\alpha} A_{\alpha}$ is a directed union of super Hopf subalgebras $A_{\alpha}$.
If the induced quotients $A_{\alpha} \to D_{\alpha}$ are all faithfully coflat, then $A \to D$ is as well.
\end{lemma}

\begin{proof}
Set $B = A^{\mathrm{co}D},\ B_{\alpha} = A^{\mathrm{co}D_{\alpha}}$; see \eqref{(3.4)}. 
Note that $B = \bigcup_{\alpha} B_{\alpha},\ D = \bigcup_{\alpha} D_{\alpha}$, directed unions. 
If $A_{\alpha} \to D_{\alpha}$ are all faithfully coflat, it follows by Proposition \ref{3.2} 
that $B_{\alpha} \hookrightarrow A_{\alpha}$ are all faithfully flat, and 
$A_{\alpha}/B_{\alpha}^{+}A_{\alpha} = D_{\alpha}$.
This implies that  
$B = \bigcup_{\alpha} B_{\alpha} \hookrightarrow \bigcup_{\alpha} A_{\alpha} = A$ is faithfully flat, and
$A/B^{+}A = \underrightarrow{\mathrm{lim}}_{\alpha}~k \otimes_{B_{\alpha}} A_{\alpha} 
= \underrightarrow{\mathrm{lim}}_{\alpha}~D_{\alpha} = D.$ 
Again by Proposition \ref{3.2}, $A \to D$ is faithfully coflat.
\end{proof}

In what follows we prove Theorem \ref{1.4}. 
Let $A \to D$ be a quotient of a super commutative Hopf algebra $A$. 
Note that $A$ is presented as a directed union $A = \bigcup_{\alpha} A_{\alpha}$ of finitely generated 
super Hopf subalgebras $A_{\alpha}$. 
Let $A_{\alpha} \to D_{\alpha}$ be the induced quotients. 
As was seen in \eqref{(4.7)}, each $\overline{D}_{\alpha}$ is a Hopf subalgebra of $\overline{D}$. 
Therefore, if $\overline{D}$ satisfies Condition (a) or (b) in the theorem, each 
$\overline{D}_{\alpha}$ is
\begin{equation}\label{(5.1)}
\mbox{(a)~cosimisimple\ \ or\ \ (b)~finite-dimensional,}
\end{equation}
respectively. In virtue of Lemma \ref{5.1}, replacing $A \to D$ with $A_{\alpha} \to D_{\alpha}$, 
we may suppose that $A$ is finitely generated, and $\overline{D}$ is in each case of \eqref{(5.1)}.

\subsection*{Proof in Case~(b).} 
By bozonization we have a surjection 
$\widetilde{A} := \mathbb{Z}_2 \cmddotltimes A \to \widetilde{D} := \mathbb{Z}_2 \cmddotltimes D$ 
of ordinary Hopf algebras. 
Since $A$ is assumed to be finitely generated and we are in case (b), Theorem \ref{4.1} implies
that $A$ and hence $\widetilde{A}$ have polynomial growth, and $\widetilde{D}$ is finite-dimensional. 
It follows by the main theorem of \cite{D} (see also \cite[Theorem 5.4]{Sk}) that $\widetilde{A}$ 
is cofree as a right 
$\widetilde{D}$-comodule, which implies that $A$ is an injective cogenerator (indeed, cofree) as 
a right $D$-comodule, as desired; see Proposition \ref{2.1}.\ $\square$

\subsection*{}

For our proof in Case (a), we need the following lemma, which immediately follows from Theorem \ref{4.1}.

\begin{lemma}\label{5.2}
Let $A$ be a super commutative Hopf algebra. 
Let ${}^{\mathrm{co}\overline{A}}A$ denote the super right coideal subalgebra of $A$ which consists
of all left $\overline{A}$-coinvariants; this is defined as in \eqref{(3.4)}, but on the opposite
side, and has the restricted $\varepsilon$ as counit.

(1) There is a counit-preserving isomorphism 
${}^{\mathrm{co}\overline{A}}A \overset{\simeq}{\longrightarrow} \wedge(W^A)$ of super algebras.

(2) Every quotient map $A \to D$ restricts to a counit-preserving surjection 
${}^{\mathrm{co}\overline{A}}A \to {}^{\mathrm{co}\overline{D}}D$ of super algebras.
\end{lemma}

\subsection*{Proof in Case (a).}
Set $R = {}^{\mathrm{co}\overline{A}}A,\ T = {}^{\mathrm{co}\overline{D}}D$. 
By Theorem \ref{4.1} (applied to $D$), $T \subset D$ is faithfully flat and 
$D/DT^{+} = \overline{D}$. 
By \cite[Proposition 1.1]{M}, we have a category equivalence,
\begin{equation}\label{(5.2)}
\mathbf{S}^{\overline{D}} \overset{\approx}{\longrightarrow} \mathbf{S}_{T}^{D},\ 
V \mapsto V \square_{\overline{D}} D,
\end{equation}
where $\mathbf{S}_{T}^{D}$ denotes the category of $(D, T)$-Hopf modules in $\mathbf{S}$.

Assume that we are in Case (a). 
Then, $\mathbb{Z}_2 \otimes \overline{D}$ is cosemisimple. 
Therefore, the category $\mathbf{S}^{\overline{D}}$, identified with 
$\mathbf{M}^{\mathbb{Z}_2 \otimes \overline{D}}$, is semisimple. 
By \eqref{(5.2)}, $\mathbf{S}_{T}^{D}$ is semisimple. It follows that each object $M$ in 
$\mathbf{S}_{T}^{D}$ is $D$-injective since the structure map $M \to M \otimes D$ 
splits in $\mathbf{S}_{T}^{D}$. 

Recall that we may suppose that $A$ is finitely generated, whence 
$\mathrm{dim}~W^A < \infty$. 
Then it follows by Lemma \ref{5.2}(2) that $A \to D$ restricts to a counit-preserving 
surjection $R \to T$ with nilpotent kernel, say $I$. 
Suppose $I^r = 0$ with $r > 0$. 
For each $0 \leq i < r$, $AI^i$ is naturally regarded as a $(D, R)$-Hopf 
module, or in notation, $AI^i \in \mathbf{S}_{R}^{D}$. 
Therefore, $AI^i/AI^{i+1} \in \mathbf{S}_{T}^{D}$, whence it is $D$-injective. 
This implies that the short exact sequence
\begin{equation}\label{(5.3)}
0 \to AI^i/AI^{i+1} \to A/AI^{i+1} \to A/AI^{i} \to 0
\end{equation}
in $\mathbf{S}_{R}^{D}$ splits $D$-colinearly. It follows that $A$ decomposes so as
\begin{equation}\label{(5.4)}
A \simeq \bigoplus_{i=0}^{r-1} AI^i/AI^{i+1}
\end{equation}
into direct summands of injective right $D$-comodules. 
Since the surjection $A/AI \to D$ induced from $A \to D$ splits in $\mathbf{S}_{T}^{D}$,
$A/AI$ is a cogenerator as well in $\mathbf{M}^{D}$, whence $A$ is an injective 
cogenerator in $\mathbf{M}^{D}$, as desired.\ $\square$

\section{Notes added after the first submission} 

\subsection{Afterthought.} 
After I posted to the arXiv the original version of this paper on January 27, 2010, 
I found that Theorems \ref{1.1} and \ref{1.4} can be generalized to a large extent, as follows,
by slightly modifying the original proofs.

Let $G, H$ be as in Theorem \ref{1.1}. Let $\overline{G}, \overline{H}$ denote the affine $k$-groups
associated to $G, H$, respectively. 
Recall that $\overline{H}$ is a affine closed subgroup of $\overline{G}$.

\begin{theorem}
The $k$-dur sheaf $G\tilde{\tilde{/}} H$ of right cosets is affine if and only if 
$\overline{G}\tilde{\tilde{/}} \overline{H}$ is affine.
\end{theorem}

In virtue of the affineness criteria for quotient dur sheaves due to Takeuchi \cite[Theorem 10]{T2} 
and Zubkov \cite[Theorem 5.2]{Z}, the theorem just formulated is translated into Hopf-algebra language 
so as (a) $\Leftrightarrow$ (c) of the next theorem.
 
Let $A \to D$ be as in Theorem \ref{1.4}. Set $\overline{A} = A_0/A_1^2$, $\overline{D} = D_0/D_1^2$; see 
\eqref{(1.2)}. 
Then there arises a natural surjection $\overline{A} \to \overline{D}$ of commutative Hopf algebras.
 
\begin{theorem}
The following are equivalent:

(a) $A \to D$ is faithfully coflat;

(b) $A \to D$ is coflat;

(c) $\overline{A} \to \overline{D}$ is faithfully coflat;

(d) $\overline{A} \to \overline{D}$ is coflat;
\end{theorem}

\begin{proof}
(a) $\Rightarrow$ (b). Trivial.

(b) $\Rightarrow$ (d). This follows since by Theorem \ref{4.1}, $A \to \overline{A}$ is faithfully coflat, 
and $D \to \overline{D}$ is coflat.

(c) $\Leftrightarrow$ (d). This is due to Doi \cite[Remark on Page 247]{Doi}.

(c) $\Rightarrow$ (a). Assume (c). As is seen from the proof of Theorem \ref{1.4} in Case (a), it 
suffices to prove the following two:
\begin{itemize}

\item[(1)]
For each $0 \leq i < r$, the $AI^i/AI^{i+1}$ in \eqref{(5.4)} is an injective object in $\mathbf{S}_{T}^{D}$. 

\item[(2)]
The surjection $A/AI \to D$ given in two lines below \eqref{(5.4)} splits in $\mathbf{S}_{T}^{D}$.
\end{itemize}

Notice from \cite[Proposition 1.1]{M} that through the category equivalence \eqref{(5.2)}, 
the last surjection $A/AI \to D$ corresponds to 
$\overline{A} \to \overline{D}$ in $\mathbf{S}^{\overline{D}}$, which splits by (c). 
This proves (2).

To prove (1), fix $i$, and set $M = AI^i/AI^{i+1}$. 
Since $A$ is an algebra in $\mathbf{S}^{D}$, we have the category $\mathbf{S}_{A}^{D}$ of $(D, A)$-Hopf 
modules in $\mathbf{S}$. Similarly we have $\mathbf{S}_{\overline{A}}^{\overline{D}}$. 
We see that $M$ has a natural structure in $\mathbf{S}_{A}^{D}$; it induces the structure of $M$ in
$\mathbf{S}_{T}^{D}$ which was used in the proof of Theorem \ref{1.4} in Case (a). 
Notice again that through \eqref{(5.2)}, $M$ corresponds to 
\begin{equation*}
M/MT^{+} = M/MR^{+}~(= AI^{i}/AI^{i}R^{+})
\end{equation*}
in\ $\mathbf{S}^{\overline{D}}$. This last object indeed has the structure in $\mathbf{S}_{\overline{A}}^{\overline{D}}$ 
which arises from the structure of $M$ in $\mathbf{S}_{A}^{D}$. 
Since by (c), there exists a unit-preserving right $\overline{D}$-colinear map $\overline{D} \to \overline{A}$, 
it follows by Doi's Theorem \cite[Theorem 1]{Doi} that every ordinary $(\overline{D}, \overline{A})$-Hopf
module is injective as a right $\overline{D}$-comodule. It follows that $M/MT^{+}$ is an injective object in 
$\mathbf{S}^{\overline{D}}$; this proves (1).
\end{proof}

\subsection{Erratum.}
In Section \ref{Introduction} the term `reductive' is confused with `linearly reductive', as was pointed out by Zubkov.
Condition (a) in Theorem \ref{1.1} and Corollary \ref{1.2} should be replaced by (a) \emph{linearly reductive}. 
Therefore, Corollary 1.2 in Case (a) does not generalize Zubkov's result cited just below the corollary.

\section*{Acknowledgements}

I would gratefully acknowledge that the present work was supported by
Grant-in-Aid for Scientific Research (C) 20540036 and Grant-in-Aid for Foreign
JSPS Fellow both from Japan Society of the Promotion of Science. 
I thank Alexandr Zubkov very much for his helpful comments.

\end{document}